\documentclass[a4paper,12pt]{article}

\usepackage{amssymb,amsmath,amsthm,latexsym,amsfonts,epsfig,euscript}

\newtheorem{fed}{Definition}[section]
\newtheorem{teo}[fed]{Theorem}
\newtheorem{lem}[fed]{Lemma}
\newtheorem{cor}[fed]{Corollary}
\newtheorem{pro}[fed]{Proposition}
\theoremstyle{definition}
\newtheorem{rem}[fed]{Remark}

\def\beq{\begin{equation}}
\def\eeq{\end{equation}}

\def\cH{\mathcal{H}}
\def\cK{\mathcal{K}}

\def\cW{\mathcal{W}}

\def\bdem{\begin{proof}}
\def\edem{\end{proof}}

\def\orto{^\perp}

\def\*A{\#\sb A}

\def\CC{\mathbb {C}}

\def\H{{\cal H}}
\def\K{{\cal K}}

\def\cH{{\cal H}}

\def\cS{{\cal S}}

\def\cM{{\cal M}}
\def\cN{{\cal N}}

\def\api{\langle}
\def\cpi{\rangle}

\def\noi{\noindent}

\def\bm{\left(\begin{array}}
\def\em{\end{array}\right)}
\def\ben{\begin{enumerate}}
\def\een{\end{enumerate}}
\def\barr{\begin{array}}
\def\earr{\end{array}}

\def\H{{\cal H}}
\def\cG{{\cal G}}
\def\lh{{L(\H)}}
\def\lh+{{L(\cH)^+}}
\def\lk+{{L(\cK)^+}}

\def\ov{\overline}

\def\noi{\noindent}

\begin{document}
\date{}


\title{{\bf Weighted generalized inverses, oblique projections and  least squares problems}
\thanks{Partially supported by CONICET (PIP 2083/00), UBACYT I030 and ANPCYT (PICT03-9521)}
}
\author{ G. Corach and A. Maestripieri
}

\maketitle

\noi {\bf Gustavo Corach (corresponding author)}

\noi Depto. de Matem\'atica, FI-UBA, \\
Paseo Col\'on 850 \\ 1063 - Buenos Aires, Argentina \\ and \\
IAM-CONICET, \\ Saavedra
15, \\ 1083  - Buenos Aires,
Argentina.\\  e-mail: gcorach@fi.uba.ar

\bigskip
\vglue 1truecm

 \noi {\bf Alejandra Maestripieri
}

\noi Instituto de Ciencias, Universidad Nacional de General
Sarmiento, \\ 1613 - Los Polvorines, Argentina \\and \\IAM-CONICET.

\noi e-mail: amaestri@ungs.edu.ar

\medskip

\vglue 1truecm

\noi {\bf Keywords:} weighted generalized inverses, oblique
projections, least squares, abstract splines.

\medskip
\noi

\vglue 1truecm
\begin{abstract}
A generalization with singular weights of Moore-Penrose generalized inverses of closed range
operators in Hilbert spaces is studied using the notion of compatibility of subspaces and positive
operators.
\end{abstract}


\section{Introduction}
Given a matrix $B \in \mathbb {C}^{m\times n}$, the Moore-Penrose
generalized inverse of $B$ is the unique matrix $C \in \mathbb
{C}^{n\times m}$ which satisfies the system
 $$ BXB=B, ~XBX=X,~
(BX)^*=BX, ~(XB)^*=XB.$$
Thus, $BC$ is the orthogonal projection
onto the column space $R(B)$ of $B$ and $CB$ is the orthogonal
projection onto the column space $R(B^*)$. In many applications,
it appears to be necessary to change the scalar products in the
spaces of input and output vectors. More precisely, given $B \in
\mathbb {C}^{m\times n}$ and $A_1 \in \mathbb {C}^{n\times n}$,
$A_2 \in \mathbb {C}^{m\times m}$ which are positive definite, the
system to be solved is

\begin{equation}\label{1}
BXB=B,~ XBX=X, ~(A_2BX)^*=A_2BX, ~(A_1XB)^*=A_1XB.\tag{$*$}
\end{equation}

Again, there exists a unique solution $C' \in
\mathbb {C}^{n\times m}$, $BC'$ (resp. $C'B$) is the orthogonal
projection onto $R(B)$ (resp. $R(B^*)$) with respect to the scalar
product on $\mathbb {C}^m$ (resp. $\mathbb {C}^n$) defined by
$A_2$ (resp. $A_1$). In some applications a singular version
of the problem needs to be solved. Thus, $A_1$ and $A_2$ are supposed to be positive
semidefinite. In this case, solutions of $(*)$ always exist but they
are infinitely many. Among them, there exists a unique solution of
minimal Euclidean norm. In other applications in which very large
numbers of variables are involved, it can be desirable to solve
system $(*)$ for bounded linear operators between Hilbert spaces.
It should be noticed that, in such cases, the first two conditions
of $(*)$ force $B$ to have a closed range and any solution will
have, also, a closed range. In this case, the existence of a
solution is not guaranteed. The goal of this paper is the complete
solution of the following problems. Let $\cH$ and $\cK$ be Hilbert
spaces, $B:\cH \to \cK$ a bounded linear operator with closed
range and $A_1: \cH \to \cH$, $A_2: \cK \to \cK$ positive
semidefinite bounded linear operators. Consider the seminorm $\|\cdot\|_{A_1}$, (resp. $\|\cdot\|_{A_2}$)
on $\H$ (resp. on $\K$) defined by  $\|x\|_{A_1}=\api A_1x,
x\cpi^{1/2}$, for $x\in \H$ ( respectively  $\|x\|_{A_2}=\api
A_2x, x\cpi^{1/2}$, for $x\in \K$).
\medskip

{\bf Problem I}

Find necessary and sufficient conditions for the existence of
solutions of system $(*)$.
\medskip

{\bf Problem II}

Find all solutions of system $(*)$, in case there exists one.

\medskip

{\bf Problem III}

Find necessary  and sufficient conditions for the existence of
$u_0\in \H$ such that $\|y-Bu_0\|_{A_2}\le \|y-Bx\|_{A_2}$ for
every $x\in\H$ and $\|u_0\|_{A_1}\le \|u\|_{A_1}$, for every $u\in
\H$ such that $\|y-Bu\|_{A_2}\le \|y-Bx\|_{A_2}$, for every
$x\in\H$.
\medskip

{\bf Problem IV}

In case there exists an $u_0$ as above, find all of them and, among them, find one of minimal Euclidean norm.
\medskip

It should be mentioned that, if the weights $A_1$ and $A_2$ are supposed to be invertible, then existence and
uniqueness of solutions of system $(*)$ follow immediately  from the analogous results on Moore-Penrose generalized
inverses, changing the inner products of $\cH$ and $\cK$ (or $\CC^n$ and $\CC^m$ in the finite dimensional case).
In this case, Problems I to IV have a unique solution. The reader is referred to the complete survey by Nashed
and Votruba \cite{[NZV]}, section 4.5, and to the more modern treatment by Nashed \cite{[N]}, with emphasis in
Banach and Hilbert space operators.
\medskip

Before the description of the main results of the paper let us give a look to the history of the subject.
\medskip

\noi {\bf Historical notes.}

The first appearance of weighted generalized inverses of matrices
is due to Greville \cite{[Gr]} who used them in problems involving
least squares fitting of curves and surfaces. As it happens with
every natural useful notion, many results on generalized inverses
have been discovered once and again by mathematicians,
statisticians and engineers. Thus, Chipman \cite{[CH1]}
reintroduced the notion for linear regression problems. Also
Goldman and Zelen \cite{[GZ]}, Watson \cite{[Wat]}, Zyskind
\cite{[Z]} and Rao and Mitra \cite{[RM]}, \cite{[RM3]},
\cite{[RM2]}, found applications to statistics. Milne \cite{[M]}
introduced a version of ``oblique pseudoinverse" for matrices and
 Ward, Boullion and Lewis \cite{[WBL]} proved that Milne's oblique pseudoinverses can be
thought as weighted generalized inverses with invertible weights.
 In a later paper \cite{[WBL2]}
 they extended some results to singular weights. In fact, in the papers mentioned above the weights are
 represented by positive definite matrices and Ward, Boullion and Lewis relaxed the hypothesis on the weights.
 Ward \cite{[War]} found a limit formula for weighted generalized inverses.
 Some related results with a different
approach have been obtained by Rao and Mitra \cite{[RM]}, \cite{[RM2]}, \cite{[MRR]} and Morley \cite{[MG-M]}.

In 1980, Eld\'en \cite{[E]} published a complete treatment of the existence of optimal weighted generalized inverses for
singular weights in finite dimensional spaces. The present paper can be seen as an extension of Eld\'en's
approach to infinite dimensional Hilbert spaces.
For recent results on this subject the reader is referred to the papers by Sun and Wei \cite{[SW]},
Stanimirovi\'c and Stankovi\'c
\cite{[SS]} and Djordjevi\'c, Stanimirovi\'c and Wei \cite{[DSW]}. For applications to parallel computing,
image processing and
many algorithmical results which use weighted generalized inverses with singular weights, the reader is referred
to the papers by Censor, Gordon and
and Gordon \cite{[CGG]}, \cite{[CGG2]} and Censor and Elfving, \cite{[CE]}, \cite{[CE2]}. The papers by Nashed and
Votruba \cite{[NZV]} and Nashed \cite{[N]}, and
the books by Rao and Mitra \cite{[RM2]} and Ben-Israel and Greville \cite{[BGT]} are excellent references,
which contain many
results on weighted generalized inverses.
\medskip

The contents of the paper are the following. Section 2 contains all results on the notion of compatibility of a
closed subspace of a Hilbert space $\cH$ and a positive bounded operator $A$ acting on $\cH$. Section 3 is devoted
to solve Problems I and II in terms of compatibility. In Section 4 we solve Problems III and IV and show an application
of our techniques by proving a result by Morley \cite{[MG-M]} on an infinite dimension regression model with
singular covariance.

\section{Preliminaries}
Throughout, $\cH,\cK,\cG$ denote Hilbert spaces, $L(\cH,\cK)$ is the space of bounded linear operators from $\cH$
to $\cK$, $L(\cH)$ is the algebra $L(\cH,\cH)$ and $L(\cH)^+$ denotes the cone of positive semidefinite  operators.
For any $C\in L(\cH,\cK)$ the image or range (resp. the nullspace) of $C$ is denoted by $R(C)$ (resp. $N(C)$). $CR(\cH,
\cK)$ is the subset of $L(\cH,\cK)$ of all operators with closed range. For any $B\in CR(\cH,\cK)$ the Moore-Penrose
inverse of $B$ is the operator $B^{\dagger} \in CR(\cK,\cH)$ such that $B^{\dagger} Bx=x$ for every $x\in N(B)^\bot$
and $B^{\dagger} y=0$ for every
$y\in R(B)\orto$. $B^{\dagger} $ is characterized by the properties $BB^{\dagger} B=B$, $B^{\dagger} BB^{\dagger}
=B^{\dagger} $, $(BB^{\dagger} )^*=BB^{\dagger} $, $(B^{\dagger} B)^*=B^{\dagger} B$.
If $\cH$ is decomposed as a direct sum of closed subspaces $\cH=\cM\oplus\cN$, the projection onto $\cM$ with nullspace
$\cN$ is denoted by $P_{\cM\|\cN}$. In particular, given a closed subspace $\cM$ of $\cH$, $P_\cM$ denotes
the projection $P_{\cM\|\cM^{\orto}}$. Denote $Q=\{Q\in L(\cH):Q^2=Q\}$. The Moore-Penrose inverse $B^{\dagger}$
of $B\in CR(\cH,\cK)$ is determined
by the properties
$BB^{\dagger} =P_{R(B)}$ and $B^{\dagger} B=P_{N(B)^{\orto}}$ or, equivalently, by:

(i) $\|BB^{\dagger} y-y\|\le \|Bx-y\|$ for every $x\in \cH$;

(ii) $\|B^{\dagger} y\|\le \|z\|$ for every $z\in \cH$ such that $\|Bz-y\|\le \|Bx-y\|$ for every $x\in\cH$.

An operator $A\in L(\cH)^+$ and a closed subspace $\cS$ of $\cH$ form a compatible
pair
$(A,\cS)$ if there exists a projection $Q\in L(\cH)$ such that $R(Q)=\cS$ and $AQ=Q^*A$. The last condition means that $Q$ is
$A$-Hermitian in the sense that $<Qx,x'>_A=$ $<x,Qx'>_A$, for every $x,x'\in \cH$, where $<x,x'>_A=<Ax,x'>$ defines a
semi inner product on $\cH$, which is an inner product only if $N(A)=\{0\}$. There is also a seminorm defined by $A$,
namely $\|x\|_A=<Ax,x>^{1/2}$ for $x\in \cH$.

Denote $P(A,\cS)=\{Q\in Q:R(Q)=\cS$ and $AQ=Q^*A\}$,  i.e., $P(A,\cS)$ is the set of $A$-Hermitian projections with fixed range
$\cS$. The set $P(A,\cS)$ can be empty (if $(A,\cS)$ is not compatible), or have one element (for example, if $A$ is positive
definite) or have infinite elements. It is easy to see that if $\cS$ is finite dimensional (and, a fortiori, if $\cH$
is finite dimensional), then every pair $A,\cS$ is compatible \cite{[CMS2]}.
The compatibility of a given pair $(A,\cS)$ has been characterized in terms of angles between subspaces and
decompositions of the ranges of $A$ and $A^{1/2}$. It has also been proven that the compatibility of $(A,\cS)$ is
equivalent to the existence of a solution of the equation $PAPX=PA(I-P)$, where $P=P_\cS$. (See \cite{[CMS]},
\cite{[CMS2]} for details).
This kind of equations can be studied applying Douglas theorem:

\begin{teo}\label{douglas}
Given
Hilbert spaces $\cH$, $\cK$, $\cG$ and operators $A\in L(\cH,\cG)$, $B\in L(\cK,\cG)$
then the following conditions are equivalent:

i) the equation $AX=B$ has a solution in $L(\cK,\cH)$;

ii) $R(B)\subseteq R(A)$;

iii) there exists $\lambda >0$ such that $BB^*\le \lambda AA^*$. In this case,
there exists a unique $D\in L(\cK,\cH)$ such that $AD=B$ and $R(D)\subseteq \ov{R(A^*)}$;
moreover, $\|D\|^2=\inf \{\lambda >0:~ BB^*\le \lambda AA^*\}$. We shall say that $D$ is the {\bf
reduced solution} of $AX=B$.
\end{teo}

The reader is referred to \cite{[D]} and \cite{[FW]} for the proof of Douglas
theorem and related results.
\medskip

Suppose that $(A,\cS)$ is compatible and consider the reduced solution $D$ of the equation $PAPX=PA(I-P)$. Define
$P_{A,\cS}=P+D$, or, in terms of the matrix representation induced by $P$,
$P_{A,\cS}=\left(\begin{matrix} 1 &d\\ 0 & 0\end{matrix}\right)$, where $P$ is identified with the identity in $L(\cS)$,
$D$
with the operator $d=D|_{\cS^\bot}\in L(\cS^\bot,\cS)$. The next theorem characterizes the set $P(A,\cS)$:

\begin{teo}
Let $A\in L(\cH)^+$ and $\cS$ a closed subspace of $\cH$ such that $(A,\cS)$ is compatible. Then $P_{A,\cS}\in P(A,\cS)$
and it
is the projection onto $\cS$ with nullspace $A^{-1}(\cS^\bot)\ominus (N(A)\cap \cS)$. The set $P(A,\cS)$ is an affine
manifold and
it can be parametrized as
$$
P(A,\cS)=P_{A,\cS} + L(\cS^\bot, N(A)\cap \cS),
$$
where $L(\cS^\bot, N(A)\cap \cS)$ is viewed as a subspace of $L(\cH)$.
\end{teo}

Given $T\in L(\cH,\cK)$, a closed subspace $\cS$ of $\cH$,
 and an element $y\in \cH$, an
{\it abstract spline} or a $(T,\cS)$-{\it spline interpolant} to $y$ is any element of the set
$$
spl(T,\cS,y)=\{x\in y+\cS:\|Tx\|\le \|T(y+s)\| \text{ for all } s\in \cS\}.
$$
It holds that $spl(T,\cS,y)=(y+\cS)\cap A(\cS)^\bot$ where
$A=T^*T$. The abstract theory of splines is due to Atteia \cite
{[Att]}. The reader is referred to \cite {[Del]} and \cite{[Gro]}
for some relationships between abstract splines and generalized
inverses.

The following theorem relates the existence of splines to compatibility:

\begin{teo}\label{2.3}
Let $T\in L(\cH,\cK)$ and $\cS$ a closed subspace of $\cH$. If $A=T^*T$, then:

a) $spl(T,\cS,y)$ is not empty for every $y\in \cH$ if and only if the pair $(A,\cS)$ is compatible.

b) If $(A,\cS)$ is compatible and $y\in \cH\setminus\cS$ then $spl(T,\cS,y)=\{(I-Q)y:Q\in P(A,\cS)\}$. Furthermore, $(I-P_{A,\cS})y$
is the unique vector in $spl(T,\cS,y)$ with minimal norm.
\end{teo}

See \cite{[CMS3]} for the proofs of these assertions.

\section{Weighted generalized inverses}
\begin{teo}Given $B \in CR(\cH,\cK)$, $A_1 \in \lh+ $ and $A_2 \in \lk+ $ there
exists $C\in L(\K,\H)$ such that

\begin{equation}\label{1}
BCB=B, CBC=C, A_1CB=B^*C^*A_1, A_2BC=C^*B^*A_2
\end{equation}
if and only if $(A_1,N(B))$ and $(A_2,R(B))$ are
compatible pairs.
\end{teo}

\bdem Suppose that $C\in L(\K,\H)$ satisfies (\ref{1}). Notice
that $C$ has closed range: in fact the
projection $P=BC$ on $\cK$ has the same range as $B$ and the
projection $Q=CB$ has the same range as $C$; of course, $Q$ is a bounded linear projection and, therefore,
its range is closed. It follows easily
that $Q$ (resp. $P$) and $B$ (resp. $C$) have the same nullspace.
Observe also that the third and fourth conditions of (\ref{1}) say
that $Q$ is $A_1-$Hermitian and $P$ is $A_2-$Hermitian. Then $I -Q$ is also $A_1$-Hermitian and
$R(I -Q)=N(Q)=N(B)$, which proves that $I -Q\in P(A_1,N(B))$. Analogously, $P\in P(A_2,R(B))$.
This shows that $(A_1,N(B))$ and $(A_2, R(B))$ are compatible pairs.

Conversely, suppose there exist $Q'\in P(A_1,N(B))$ and $P\in P(A_2, R(B))$. Then $Q=I -Q'$ is $A_1$-Hermitian
and $N(Q)=R(Q')=N(B)$.
Consider the decomposition $\cK=R(B)\oplus N(P)$ and define $C:\cK\to \cH$ by $C(Bx+z)=Qx$, for $x\in \cH$, $z\in
N(P)$.
$C$ is well defined because $N(B)=N(Q)$. It is also easy to check that $C$ is a linear operator, with $R(C)=R(Q)$ and
$N(C)=N(P)$; $C$ is also bounded, because $B|_{R(Q)}:R(Q)\to R(B)$ is an isomorphism by the closed graph theorem
and $C|_{R(B)}=(B|_{R(Q)})^{-1}$. This also implies $BCB=B$. It remains to prove the other conditions of (\ref{1}).
On one side, it holds
$CBC(Bx+z)=CBx=C(Bx+z)$ for every $x\in\cH$ and $z\in N(P)$. On the other side, $CB=Q$ is $A_1$-Hermitian and $BC=P$ is
$A_2$-Hermitian.
\edem

\medskip
From now on, $GI(B,A_1,A_2)$ denotes the set of all bounded linear solutions of $(*)$:
$$
GI(B,A_1,A_2)=\{C\in CR(\cK,\cH):BCB=B, CBC=C,A_1CB=B^*C^*A_1,A_2BC=C^*B^*A_2\}.
$$
The proof of the theorem above and the characterization of the set of generalized Hermitian projections of a given
range described in section 2, provide the following parametrization of $GI(B,A_1,A_2)$:

\begin{pro}
The set $GI(B,A_1,A_2)$ is parametrized by the vector space
$$
L(N(B)^{\orto}, N(A_1)\cap N(B))\times L(R(B)^{\orto}, N(A_2)\cap R(B)).
$$
\end{pro}

\bdem
The proof of the theorem above shows that the construction of a bounded linear solution of ($*$), if there
exists any, is based in the choice of two projections, namely, $I-Q\in P(A_1,N(B))$ and $P\in P(A_2,R(B))$.
It is not difficult to prove that different choices provide different solutions of ($*$). On the other hand,
following the notations and results of section 2, $P(A_1, N(B))$ is in bijection with $L(N(B)^{\orto}, N(A_1)\cap
N(B))$ and $P(A_2, R(B))$ is in bijection with $L(R(B)^{\orto}, N(A_2)\cap R(B))$. With these comments, the result
follows straightforward.
\edem

The parametrization just obtained is quite indirect. The following results of this section are devoted to find
more explicit parametrizations of $GI(B, A_1,A_2)$.

The first goal is to generalize Douglas theorem in order to get convenient solutions of Douglas-type equations.

\begin{teo} Let $\cH,\cK$ and $\cG$ be Hilbert spaces. Given $A\in L(\cH,\cG)$ and $B\in L(\cK,\cG)$ such that
$R(B)\subseteq R(A)$, for every closed subspace $\cM$ of $\cH$ such that $\cH=N(A)\oplus \cM$ there exists a unique
solution
$C\in L(\cK,\cH)$ of the operator equation $AX=B$ such that $R(C)\subseteq \cM$. The nullspace of $C$ coincides with
that of $B$.
\end{teo}

\bdem
Consider the reduced solution $C'\in L(\cK,\cH)$ of $AX=B$ and define $C=P_{\cM\|N(A)}C'$. Obviously, $R(C)\!
\subseteq\!\cM$. Observe that $AC=AP_{\cM\|N(A)}C'\!=\!AC'=\!B$ because  $A(I-P_{\cM\|N(A)})=AP_{N(A)\|\cM}=0$.
Therefore, $AC=B$, which proves the existence statement.

Suppose that $D\in L(\cK,\cH)$ satisfies $AD=B$ and $R(D)\subseteq\cM$. Then, $A(D-C)=0$ so that $R(D-C)\subseteq N(A)$.
But $R(D-C)\subseteq \cM$ and, therefore, $R(D-C)\subseteq N(A)\cap \cM=\{0\}$. This shows that $D=C$.

The last assertion follows easily: $N(C)\subseteq N(B)$ because $AC=B$; conversely, if $Bx=0$ then $Cx\in N(A)\cap
R(C)\subseteq N(A)\cap\cM=\{0\}$, which shows that $N(B)\subseteq N(C)$.
\edem

\medskip
Given $B\in CR(\cH,\cK)$ let us denote $B\{1\}=\{C\in L(\cK,\cH):BCB=B\}$ and $B\{1,2\}=\{C\in L(\cK,\cH):BCB=B$ and
$CBC=C\}$.
Following the notations of Ben Israel and Greville \cite{[BGT]}, we call any $C\in B\{1\}$ an $\{1\}$-inverse of $B$ and
any $C\in B\{1,2\}$ an $\{1,2\}$-inverse of $B$.

\begin{cor}\label{3.4}
Consider $B\in CR(\cH,\cK)$ and projections $Q\in L(\cH)$, $P\in L(\cK)$ such that $N(Q)=N(B)$ and $R(P)=R(B)$. Then
there exists a unique solution $C\in L(\cK,\cH)$ of
\begin{equation}\label{2}
BX=P, ~ R(X)=R(Q).
\end{equation}
It holds $C\in B\{1,2\}$ and $N(C)=N(P)$.
\end{cor}

\bdem
Observe the decompositions $\cH=N(Q)\oplus R(Q)=N(B)\oplus R(Q)$ and the inclusion $R(P)\subseteq R(B)$.
By Theorem 3.3, there exists a unique $C\in L(\cH,\cK)$ such that $BC=P$ and $R(C)\subseteq R(Q)$, and $C$ satisfies
also $N(C)=N(P)$. It remains to prove that $C\in B\{1,2\}$ and $R(Q)\subseteq R(C)$.

Since $R(P)=R(B)$ it follows that $PB=B$ so that $BCB=PB=B$; also $C(I-P)=0$ because $N(C)=N(P)=R(I-P)$; therefore,
$CBC=CP=C$ and this proves that $C\in B\{1,2\}$. In order to prove the inclusion $R(Q)\subseteq R(C)$, observe first
that $N(B)\subseteq N(CB)\subseteq N(CBC)=N(B)$, so that $N(CB)=N(B)=N(Q)$. Then, $CB$ and $Q$ are bounded
linear projections with the same nullspace and $R(CB)\subseteq R(C)\subseteq R(Q)$ and, therefore, $CB=Q$ and,
a fortiori, $R(C)=R(Q)$.
\edem

Observe first that any solution $C$ of
\begin{equation}\label{2'}
BX=P, ~N(X)=N(P)\tag{$2'$}
\end{equation}
satisfies $BCB=B$ because $BC=P$ and $PB=B$; similarly,
$C(I-P)=0$
because $N(C)=N(P)$ and then, $CBC=CP=C$. Thus $C\in B\{1,2\}$. By the generalization of Douglas theorem, there exists a
unique solution $C\in L(\cK,\cH)$  of
\begin{equation}\label{3}
BX=P, ~N(X)=N(P), ~R(X)\subseteq R(Q).
\end{equation}

By the first remark, it holds $CBC=C$, so that $(CB)^2=CB$. Moreover, $N(C)=N(BC)=N(P)$ and $R(CB)\subseteq R(C)
\subseteq R(Q)$. Thus, $CB$ is a projection with $N(CB)=N(B)=N(Q)$ and $R(B)\subseteq R(Q)$. By the first remark, it
holds $CBC=C$, so that $(CB)^2=CB$. Moreover, $N(C)=N(BC)=N(P)$ and $R(CB)\subseteq R(C)\subseteq R(Q)$. Thus, $CB$
is a projection with $N(CB)=N(B)=N(Q)$ and $R(B)\subseteq R(Q)$. By elementary theory of projections, it holds $CB=Q$
and, a fortiori, $R(Q)=R(CB)\subseteq R(C)\subseteq R(Q)$, so that $R(C)=R(Q)$ and this proves that $C$ is a solution
of (\ref{2}). Uniqueness of solutions of (\ref{2}) follows from that of (\ref{3}).

Notation: In what follows, $B^{\dagger} _{P,Q}$ denotes the unique solution of (\ref{2}). It follows from the proof that
$B^{\dagger} _{P,Q}$ is the unique operator in $L(\cK,\cH)$ such that
\begin{equation}\label{4}
BB^{\dagger} _{P,Q}=P \text{ and } B^{\dagger} _{P,Q}B=Q.
\end{equation}

Of course, $B^{\dagger} _{P,Q}$ has closed range, namely $R(Q)$.

As a corollary, we get another parametrization of $GI(B,A_1,A_2)$:

\begin{cor}\label{3.5}
Suppose that $(A_1,N(B))$ and $(A_2, R(B))$ are compatible pairs. Then
$$
GI(B,A_1,A_2)=\{B^{\dagger} _{P,I-Q}:P\in P(A_2,R(B)), Q\in P(A_1, N(B))\}.
$$
\end{cor}

\bdem
Let $Q\in P(A_1,N(B))$, $P\in P(A_2,R(B))$. Then $B^{\dagger} _{P,I-Q}$ satisfies the equivalent of (\ref{4}):
$BB^{\dagger} _{P,I-Q}=P$,
$B^{\dagger} _{P,I-Q}B=I -Q$. The fact that $P$ (resp. $I -Q$) is $A_2$ (resp. $A_1$)-Hermitian, together with
the identities $BB^{\dagger} _{P,I-Q}B=B$ and $B^{\dagger} _{P,I-Q}BB^{\dagger} _{P,I-Q}=B^{\dagger}_{P,I-Q}$,
prove that $B^{\dagger}_{P,I-Q}$ belongs to $GI(B,A_1,A_2)$.

Conversely, if $C\in GI(B,A_1,A_2)$ then $C$ satisfies ($*$). Then $P=BC\in P(A_2,R(B))$ and if $Q=CB$ then $I
-Q\in P(A_1,N(B))$. Thus, $C$ is the unique solution $B^{\dagger} _{P,I-Q}$ of (\ref{2}).
\edem

The next result gives a better way of constructing $B^{\dagger} _{P,Q}$ in terms of $B\{1\}$. As a corollary we shall
get a simpler parametrization of $GI(B, A_1,A_2)$.

\def\B1{B^{(1)}}

\begin{pro} Given $B\in CR(\cH,\cK)$ and projections $Q\in L(\cH)$ and $P\in L(\cK)$ such that $N(Q)=N(B)$ and $R(P)
=R(B)$ it holds $B^{\dagger} _{P,Q}=Q\B1 P$ for any $\B1\in B\{1\}$.
\end{pro}

\bdem
Take any $\B1\in B\{1\}$ and let $\cM=R(\B1 B)$. Then $\cH=\cM\oplus N(B)$ because $N(B)=N(\B1 B)$, and $\B1 B$ is a
projection onto $\cM$. Define $C=Q\B1 P$. Straightforward computations show that $N(C)=N(P)$. Let us prove that $R(C)=
R(Q)$: observe that $R(Q)=Q\cM$; then $R(C)=Q(R(B'P))=Q(R(B'B))=Q\cM$, because $P$ and $B$ have the same range.
Finally, the identity $BQ=B$, due to the fact that $N(Q)=N(B)$, implies $BC=BQ\B1 P=B\B1 P=P$, because $B\B1$ is a
projection onto $R(B)=R(P)$. Thus, $C$ satisfies (\ref{2}) and, by Proposition \ref{3.4} it follows that
$B^{\dagger} _{P,Q}=Q\B1 P$, as claimed.
\edem

The last result of this section gives a more explicit parametrization of $GI(B,A_1,A_2)$. The fact that we use the
Moore-Penrose inverse of $B$ instead of an arbitrary choice of a $\{1\}$-inverse of $B$ is not relevant.

\begin{teo}
If $(A_{1},N(B))$ and $(A_2, R(B))$ are compatible pairs then
$$
GI(B,A_1,A_2)=\{(I-Q)B^{\dagger}P:Q\in P(A_1,N(B)), P\in P(A_2,R(B))\}
$$
where $B^{\dagger}$ is the Moore-Penrose inverse of $B$.
\end{teo}

\bdem
It follows by combining the last proposition with Corollary \ref{3.5}.
\edem

\def\seq{\subseteq}
\def\cK{{\cal K}}
\def\cL{{\cal L}}
\def\cV{{\cal V}}
\def\cW{{\cal W}}

\begin{rem}
Milne \cite{[M]} defined what he called the {\it oblique pseudoinverse} of an operator acting between finite dimensional
Hilbert spaces. Let $B\in L(\cV,\cW)$ and let $\cK\seq \cV$, $\cL \seq \cW$ be two subspaces such that $\cV=\cK
\oplus N(B)$ and
$\cW=R(B)\oplus \cL$. The oblique pseudoinverse of $B$ with respect to the subspaces $\cK$ and $\cL$ is defined as the
unique $B^{\dagger}_{\cK,\cL}\in L(\cW,\cV)$ satisfying $B^{\dagger}_{\cK,\cL}Bv=v$ for every $v\in \cK$ and
$B^{\dagger}_{\cK,\cL}w=0$
for every $w\in \cL$. Milne's definition and results have trivial extensions to closed range operators between infinite
dimensional Hilbert spaces. If $Q\in L(\cV)$ is the projection onto $\cK$ with nullspace $N(B)$ and $P\in L(\cW)$ is the
projection onto $R(B)$ with nullspace $\cL$, then it can easily be shown that $B^{\dagger}_{\cK,\cL}$ satisfies
(\ref{2}) so
that $B^{\dagger}_{\cK,\cL}=B^{\dagger}_{P,Q}$. It should be remarked that Milne proved that
$B^{\dagger}_{\cK,\cL}=Q\B1 P$ for any $\B1\in B\{1\}$. An algebraic treatment of the properties of $B^\dagger_{P,Q}$
can be found in the survey by Nashed and Votruba \cite{[NZV]}.
\end{rem}

\section{Least squares formulation}

The great impact that Moore-Penrose inverses have in science is due to the fact that they solve a least squares
problems, namely, $B^\dagger c$ is the unique vector in $\cH$ with minimal norm among those which minimize $\|Bx-c\|$.
We generalize this result for the weighted case, i.e., if we consider weights $A_1$ and $A_2$ on $\cH$ and $\cK$,
respectively.

\begin{fed}
Given $B\in L(\cH,\cK)$, $A_1\in L(\cH)^+$, $A_2\in L(\cK)^+$ and $y\in\cK$, an element $u\in\cH$ is said
to be an $A_2$-least squares solution (hereafter, $A_2-LSS$) of the equation
\begin{equation}\label{5}
Bx=y
\end{equation}
if $\|Bu-y\|_{A_2}\le \|Bx-y\|_{A_2}$
for every $x\in\cH$.

An element $u_0\in \cH$ is said to be an $A_1A_2$-least squares solution $(hereafter A_1A_2-LSS)$ of (\ref{5}) if
$u_0$ is an $A_2$-LSS of (\ref{5}) and $\|u_0\|_{A_1}\le \|u\|_{A_1}$ for every $u$ which is an $A_2$-LSS of (\ref{5}).
\end{fed}

\begin{lem} Given $B \in CR(\cH,\cK)$ and $ A \in \lk+ $  there exists
an $A-$LSS $u \in \H$ of the equation $Bx=y$ for every $y\in\K$ if
and only if the pair $(A, R(B))$ is compatible.
\end{lem}

\bdem Observe that $u\in \H$ is an  $A-$LSS  of $Bx=y$ if and only
if $u\in spl(A^{1/2},R(B), y)$. Then $Bx=y$ admits an $A-$LSS for every $y\in\cK$ if
and only if $spl(A^{1/2},R(B), y)$ is not empty, for every $y\in
\cK$, which, by item a) of Theorem \ref{2.3}, is equivalent to the compability of
$(A,\cal S)$.
\edem

\begin{rem}
Given $B\in CR(\cH,\cK)$, $T\in L(\cK,\cG)$ and $y\in \cK$, it follows from the preliminaries that $u_0\in\cH$
is an $A$-LSS of $Bx=y$ (where $A=T^*T\in L(\cK)^+)$ if and only if $y-Bu_0\in spl (T,R(B),y)$. Therefore, by
the characterization of splines in the preliminaries section, given $y\in \cK$ there exists an $A$-LSS $u_0$ of
$Bx=y$ if and only if $y\in Bu_0+R(AB)^{\orto}$.

The next result determines all $A$-LSS of (\ref{5}) if $(A,R(B))$ is compatible.
\end{rem}

\begin{pro}\label{4.4}
Given $B\in CR(\cH,\cK)$, $A\in L(\cK)^+$ such that $(A,R(B))$ is compatible, and $y\in \cK\setminus R(B)$, then
$u\in \cH$ is an $A$-LSS of (\ref{5}) if and only if there exists $P\in P(A,R(B))$ such that $Bu=Py$.
\end{pro}

\bdem
If $y\in \cK\setminus R(B)$ then by Theorem \ref{2.3} it holds $spl(A^{1/2},R(B),y)=\{(I-P)y:P\in P(A,R(B))\}$ so,
by the last
remark, $u$ is an $A$-LSS of (\ref{5}) if and only if $y-Bu\in spl(A^{1/2},R(B),y)$ and the result follows.
\edem

Observe that if $y\in R(B)$ then every element $u\in B^{-1}(\{y\})$ (i.e., every solution of $Bx=y$) is trivially
an $A$-LSS solution of (\ref{5}). In fact, in this case, $u\in\cH$ is an $A$-LSS of (\ref{5}) if and only if $y-Bu\in N(A)$.

\begin{rem}\label{4.5}
If $N(A)\cap R(B)=\{0\}$ then $u\in\cH$ is an $A$-LSS of (\ref{5}) if and only if $B^*A(Bx-y)=0$. In fact,
$P(A,R(B))$ consists of a single element $P=P_{A,R(B)}$ whose nullspace is $A^{-1}(R(B)^{\orto})=R(AB)^{\orto}$.
Straightforward computations prove the statement.
\end{rem}

\medskip

\def\setm{\setminus}

If $y\in\cK\setm R(B)$ then, by Proposition \ref{4.4}, the set of all $A$-LSS of (\ref{5}) is given by $\bigcup\{B^{-1}
\{Py\}:P\in P(A, R(B))\}$ and, for a fixed $P$, $B^{-1}\{Py\}=x_0+N(B)$, where $x_0$ is the unique element of
$B^{-1}\{Py\}\cap N(B)^{\orto}$. Notice that $x_0=P_{N(B)^\bot}x_0=B^\dagger Bx_0=B^\dagger Py$.

\medskip
Let us study a minimizing problem in $\cH$.

\begin{lem}\label{4.6}
Consider $B\in L(\cH,\cK)$ and $A\in L(\cH)^+$ such that $(A,N(B))$ is compatible.
Then, for every non zero $x_0\in N(B)^\bot$
and $u\in x_0+N(B)$, it holds $\|u\|_A\le\|x\|_A$ for every $x\in x_0+N(B)$ if and only if there exists $Q\in P(A,N(B))$
such that $u=(I-Q)x_0$.
\end{lem}

\bdem
Decompose $u=x_0+P_{N(B)}u$. Then $\|u\|_A\le \|x\|_A$ for every $x\in x_0+N(B)$ if and only if $\|u\|_A\le
\|x_0+P_{N(B)}x\|_A$ for every $x\in \cH$ or equivalently, $u$ is an $A$-LSS of the equation $P_{N(B)}x=-x_0$.
Applying the last proposition to the operator $P_{N(B)}$ and
the vector $x_0\in\cH\setm N(B)=\cH\setm R(P_{N(B)})$ this is equivalent to the existence of
$Q\in P(A, N(B))$ such that $P_{N(B)}u=-Qx_0$. Adding $x_0$ to the last equality, we get $u=(I-Q)x_0$ as claimed.
\edem

\begin{rem}
If $x_0=0$ then $\|u\|_A\le \|x\|_A$ for every $x\in N(B)$ if and only if $\|u\|_A=0$, which means that $u\in N(A)\cap
N(B)$.
\end{rem}

\medskip

We are now in position of finding all $A_1A_2$-LSS of $Bx=y$.

\begin{pro}\label{4.8}
Let $B\in CR(\cH,\cK)$, $A_1\in L(\cH)^+$ and $A_2\in L(\cK)^+$ be such that $(A_1, N(B))$ and $(A_2,R(B))$ are compatible
pairs. Consider $y\in \cK\setm R(B)$ and $u\in \cH\setm N(B)$. Then $u$ is an $A_1A_2$-LSS of the equation $Bx=y$ if
and only if there exist $Q\in P(A_1,N(B))$ and $P\in P(A_2,R(B))$ such that $u=(I-Q)B^\dagger Py$.
\end{pro}

\bdem
Suppose that $u$ is an $A_1A_2$-LSS of $Bx=y$. In particular, $u$ is an $A_2$-LSS of $Bx=y$ and, by Proposition \ref{4.4}
there
exists $P\in P(A_2,R(B))$ such that $Bu=Py$. Then $x_0=P_{N(B)^\bot}u=B^\dagger Bu=B^\dagger Py$ is non zero because $Bu
\ne 0$. By the lemma above, replacing $A$ by $A_1$, there exists $Q\in P(A_1,N(B))$ such that $u=(I-Q)x_0=(I-Q)B^\dagger
Py$. Conversely, suppose $u=(I-Q)B^\dagger Py$ for some $Q\in P(A_1,N(B))$ and $P\in P(A_2, R(B))$.
Then $Bu=B(I-Q)B^\dagger
Py=BB^\dagger Py=Py$ (the second equality holds because $BQ=0$; the third one follows from the facts that $BB^\dagger
=P_{R(B)}$ and $P$ projects onto $R(B)$). Then, by Proposition \ref{4.4}, $u$ is an $A_2$-LSS of $Bx=y$. On the other hand
$u=(I-Q)B^\dagger Py=B^\dagger Py=QBPy$ is the decomposition of $u$ according to $\cH=N(B)^\bot\oplus N(B)$ and from the
lemma above it follows that $\|u\|_A\le \|z\|_{A_1}$ for every $z\in B^\dagger Py+ N(B)$, which is the set of $A_2$-LSS
of $Bx=y$, by the comments following Proposition \ref{4.4}. This finishes the proof.
\edem

\begin{teo}
Given $A_1,A_2,B$ and $y$ as before consider the problem
\begin{equation}\label{6}
\min\{\|y-Bu\|:u \text{ is an } A_2-\text{LSS of } Bx=y\}
\end{equation}
Then:

i) $u_0$ is a solution of (\ref{6}) if and only if $Bu_0=P_{A_2,R(B)}y$;

ii) $u_0$ is a solution of (\ref{6}) and an $A_1A_2$-LSS of $Bx=y$ if and only if $u_0=(I-Q)B^\dagger P_{A_2,R(B)}y$
for some $Q\in P(A_1, N(B))$;

iii) the unique minimal norm element of the set $\{(I-Q)B^\dagger P_{A_2,R(B)}y: Q\in P(A_1,N(B))\}$ is $(I-P_{A_1, N(B)})
B^\dagger P_{A_2,R(B)}y$.

\end{teo}

\bdem
To prove i) observe that by Proposition \ref{4.4} $u_0$ is an $A_2$-LSS of $Bx=y$ if and only if there exists $P\in
P(A_2,R(B))$ such that $Bu=Py$; then we look for
$$
\min\{\|(I-P)y\|:P\in P(A_2,R(B))\}.
$$
But, by theorem \ref{2.3} in the Preliminaries, this minimum is attained in $(I-P_{A_2,R(B)})y$ so that $Bu_0=
P_{A_2,R(B))}y$.

In a similar way, by proposition \ref{4.8} and i), $u_0$ is an $A_1A_2$-LSS of $Bx=y$ and a solution of (\ref{6}) if
and only if there exists $Q\in P(A_1,N(B))$ such that $u_0=(I-Q)B^\dagger P_{A_2,R(B)}y$ and ii) follows.

To prove iii) observe that the minimum of the set
$$
\{\|(I-Q)B^\dagger P_{A_2,R(B)}\|:Q\in P(A_1,N(B))\}
$$
is attained, by Theorem \ref{2.3}, in $(I-P_{A_1,N(B)})B^\dagger P_{A_2,R(B)}$.
\edem


In \cite{[MG-M]} Morley  solved the following problem: Given a (densely defined unbounded) linear operator $B:\cH\to\cK$,
with $R(B^*)$ closed, $c\in R(B^*)$ and $V\in L(\cK)$ such that $V^2$ positive semidefinite, find
\begin{equation}\label{7}
\min\{<V^2y,y>:B^*y=c\}.
\end{equation}

If $g$ is a solution of this minimizing problem, $g$ is called a {\it best linear unbiased estimator } (BLUE).

This result is equivalent to solving the following least squares problem with linear equality constraints:
given $C\in L(\cH,\cK)$,  a closed subspace $\cS$ of $\cH$, $x_0\in \cH$ and $y\in \cK$, find
$$
\inf\{\|Cx-y\| : x\in x_0+\cS\}.
$$

In fact,
$$
\min\{<V^2y,y>:B^*y=c\}=\min\{\|y\|_{V^2}:B^*y=c\}.
$$

Observe that $B^*y=c$ if and only if there exists $w\in N(B^*)$ such that $y=B^{*\dagger} c+w$; so that
(\ref{7}) is equivalent to the problem of finding
$$
\min\{\|B^{*\dagger} c+w\|_{V^2} : w\in N(B^*)\}.
$$

The next proposition is a proof of Morley's result in terms of compatible pairs, for the case of bounded operators.

\begin{pro}\label{5.1}
Consider $C\in CR(\cH,\cK)$, $\cS$ a closed subspace of $\cH$, $x_0\in\cH$ and $y\in\cK$ such that the pair $(C^*C,\cS)$
is compatible. Then $\|Cu-y\|\le \|Cx-y\|$ for every $x\in x_0+\cS$ if and only if there exists $Q\in P(C^*C,\cS)$
such that $u=(I-Q)(x_0-C^\dagger y)$.
\end{pro}

\bdem
Observe that $\|Cx-y\|^2=\|Cx-P_{R(C)}y\|^2 + \|P_{R(C)^\bot}y\|^2$ so that
$$
\inf_{x\in x_0+\cS} \|Cx-y\|^2=\|P_{R(C)^\bot}y\|^2 +\inf_{x\in x_0+\cS}\|Cx-P_{R(C)}y\|^2 .
$$

If $u_0=C^\dagger P_{R(C)}y=C^\dagger y$ then $\|Cx-P_{R(C)}y\|=\|C(x-C^\dagger y)\|=\|x-C^\dagger y\|_{C^*C}=\|x-u_0\|_{C^*C}$

\noi so that
$$
\inf_{x\in x_0+\cS}\|Cx-y\|^2=\|P_{R(C)^\bot}y\|^2 +\inf_{x\in x_0+\cS}\|x-x_0\|^2_{C^*C}=
\|P_{R(C)^\bot}y\|^2 +\inf_{x\in x_0-u_0+\cS}\|x\|^2_{C^*C}.
$$

By Lemma \ref{4.6} it follows that $\|u\|_{C^*C}\le \|x\|_{C^*C}$ for every $x\in x_0-u_0+\cS$ if and
only if there exists $Q\in P(C^*C,\cS)$ such that $u=(I-Q)(x_0-u_0)=(I-Q)(x_0-C^\dagger y)$.
\edem


\medskip

\end{document}